\documentclass[12pt]{article}    
\usepackage{latexsym,graphicx,amssymb,curves,amsmath}

\newtheorem{thm}{Theorem}[section]

\newtheorem{lem}[thm]{Lemma}

\newtheorem{defn}[thm]{Definition}

\def\proof{\noindent\textbf{Proof}.\ \ }
\def\qed{\hfill$\blacksquare$ \\}
\newcounter{egctr}

\begin{document}


\title{\bf On Local Antimagic Vertex Coloring for Corona Products of Graphs}


%
\author{S. Arumugam\dag, Yi-Chun Lee\ddag, K. Premalatha\dag, Tao-Ming Wang\ddag\footnote{Corresponding author with contact address wang@go.thu.edu.tw and the research is supported partially by the National Science
Council of Taiwan under Grant MOST 106-2115-M-029-001-MY2.} \\~~\\~~\\
\dag National Centre for Advanced Research in Discrete Mathematics\\
Kalasalingam University\\
Anand Nagar, Krishnankoil, 626 126 Tamil Nadu, India\\~~\\
\ddag Department of Applied Mathematics \\
Tunghai University\\
Taichung, 40704 Taiwan
}

\maketitle

%
%
\begin{abstract}
Let $G = (V, E)$ be a finite simple undirected graph without $K_2$ components. A
bijection $f : E \rightarrow
\{1, 2,\cdots, |E|\}$ is called a {\bf local antimagic labeling} if for any two
adjacent vertices $u$ and $v$, they have different vertex sums, i.e. $w(u) \neq w(v)$, where the vertex sum $w(u) =
\sum_{e \in E(u)} f(e)$, and $E(u)$
is the set of edges incident to $u$. Thus any local antimagic labeling induces a proper
vertex coloring of $G$ where the vertex $v$ is assigned the color(vertex sum) $w(v)$. The {\bf local antimagic
chromatic number} $\chi_{la}(G)$ is the minimum number of colors taken over all colorings
induced by local antimagic labelings of $G$.
In this article among others we determine completely the local antimagic chromatic number $\chi_{la}(G\circ \overline{K_m})$ for the corona product of a graph $G$ with the null graph $\overline{K_m}$ on $m\geq 1$ vertices, when $G$ is a path $P_n$, a cycle $C_n$, and a complete graph $K_n$.
\end{abstract}

\smallskip{}




\section{Background and Introduction}
By a graph $G = (V, E)$ we mean a finite simple undirected graph without $K_2$ components.
For graph theoretic terminology we refer to Chartrand and Lesniak
\cite{CL}. Hartsfield and Ringel\cite{Ad} introduced the concept of antimagic labeling of a graph, and conjectured that every connected graph except $K_2$ admits such an antimagic labeling, which remains unsettled till today.
\begin{defn}
Let $G = (V, E)$ be a finite simple undirected graph without isolated vertices. Let $f : E \rightarrow \{1, 2, \cdots , |E|\}$ be a bijection. For each
vertex $u \in V(G)$, the weight $w(u) = \sum_{e \in E(u)} f(e)$, where $E(u)$ is the set of edges
incident to $u$. If $w(u) \neq w(v)$ for any two distinct vertices $u$ and $v \in V(G)$, then
$f$ is called an {\bf antimagic} labeling of $G$. A graph $G$ is called antimagic if $G$ has an
antimagic labeling.
\end{defn}
Furthermore $f$ is called a {\bf local antimagic labeling} if $w(u) \neq w(v)$ only for any two {\bf adjacent} vertices $u$ and $v \in V(G)$.
Thus it is clear that any local antimagic labeling must be antimagic and induces a proper
vertex coloring of $G$ where the vertex $v$ is assigned the color(vertex sum) $w(v)$.

This notation {\bf local antimagic labeling} was raised in 2017 by the following two sets of authors independently: S. Arumugam, K. Premalatha, M. Ba\v{c}a, A. Semani\v{c}ov\'a-Fe\v{n}ov\v{c}\'ikov\'a\cite{A}, and J. Bensmail, M. Senhaji and K. S. Lyngsie\cite{B}.
Both groups also raised the following conjecture: Every connected graph other than $K_2$ is local antimagic. However this conjecture has been proved by J. Haslegrave, using the probabilistic method, in \cite{H} more recently.

The {\bf local antimagic
chromatic number} $\chi_{la}(G)$ is the minimum number of colors taken over all colorings
induced by local antimagic labelings of $G$.
Let $\chi(G)$ be the usual chromatic number of a graph $G$.
For any graph $G$, $\chi_{la}(G) \ge \chi(G)$. It was noted by S. Arumugam et al.\cite{A} that the difference $\chi_{la}(G) - \chi(G)$ can
be arbitrarily large as shown in the following theorem.
\begin{thm}\label{leaf-lower-bound}{\bf (\cite{A}, 2017) }\\
For any tree $T$ with $l$ leaves, $\chi_{la}(T) \ge l+1$.
\end{thm}
Also in \cite{A} the local antimagic chromatic number of paths, cycles, friendship graphs, wheels, complete bipartite graphs were studied.

\begin{figure}[h]
\centering
\includegraphics[width=13cm]{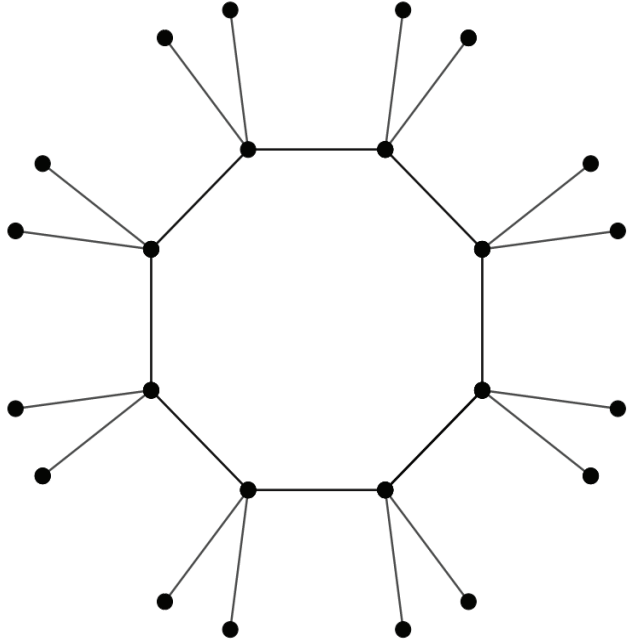}
\caption{The corona product $C_{8} \circ \overline{K_2}$}
\label{C8_K2}
\end{figure}

In this article we study the local antimagic vertex coloring of the corona product for graphs.
The {\bf corona product} of two graphs $G$ and $H$ is the graph $G \circ H$ obtained by taking one
copy of $G$ along with $|V(G)|$ copies of $H$, and via putting extra edges making the $i$-th vertex of $G$ adjacent to every vertex of the $i$-th copy of $H$, where $1 \leq i \leq |V(G)|$. See Figure~\ref{C8_K2} for example the corona product $C_{8} \circ \overline{K_2}$ for an 8-cycle $C_8$ and the null graph $\overline{K_2}$ on two vertices. We denote $\overline{K_m}$ as the complement graph of $K_m$, which is called the null graph on $m$ vertices. Also note that $\overline{K_1} = K_1$.
Note that in general the corona product operation is not commutative, i.e., $G \circ H$ is not the same with $H \circ G$. The corona product for graphs was introduced by Frucht and Harary \cite{FH} in 1970.

In the following sections, we determine completely the local antimagic chromatic number $\chi_{la}(G\circ \overline{K_m})$ for the corona product of a graph $G$ with the null graph $\overline{K_m}$, when $G$ is a path $P_n$, a cycle $C_n$, and a complete graph $K_n$.

\section{Corona Product of Paths with Null Graphs}
\subsection{The Chromatic Number of $P_{n} \circ K_1$}

\begin{figure}[h]
\centering
\includegraphics[width=14cm]{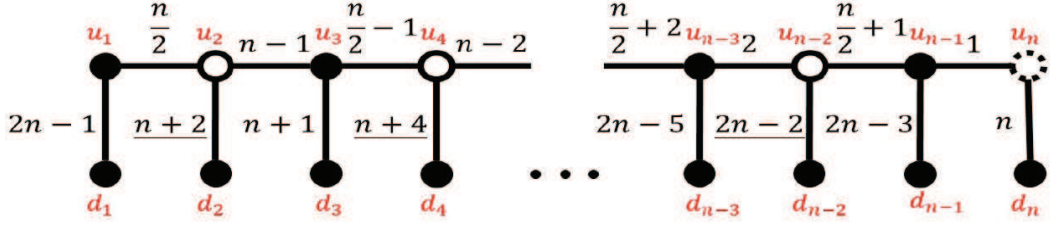}
\caption{The edge labeling of $P_{n} \circ K_1$ for $n \equiv 0~(mod~4)$}
\label{even-for-8}
\end{figure}

Let us start with the graph corona product $P_{n}$ with $K_1$ as in Figure~\ref{even-for-8}, where the vertex set is $V(P_{n} \circ K_1) = \{ u_1,u_2, \cdots, u_n, d_1,d_2, \cdots, d_n \}$ and the edge set is $E(P_{n} \circ K_1) = \{ u_i u_{i+1} : ~ 1 \leq i \leq n-1 \} \cup \{ u_i d_i : ~1 \leq i \leq n \}$. In order to get the local antimagic chromatic number $\chi_{la}(P_{n} \circ K_1)$ of $P_{n} \circ K_1$, we study the upper bound and lower bound as in the following.

\begin{lem}
$\chi_{la} ( P_{n} \circ K_1 )\leq n+2$ for $n \equiv 0~ (mod~ 4)$ and $n \ge 4$.
\end{lem}
\proof

Label the edges $u_{i} u_{i+1}, 1\leq i \leq n-1,$ as follows:

\begin{displaymath}
f(u_{i} u_{i+1})= \left\{ \begin{array}{ll}
\frac{n}{2} - \frac{i-1}{2},~~1\leq i \leq n-1, i~ \mbox{odd.}\\
n - \frac{i}{2},~~2\leq i \leq n-2, i~ \mbox{even.}
\end{array} \right.
\end{displaymath}

Label the edges $u_{i} d_{i}, 1\leq i \leq n,$ as follows:

\begin{displaymath}
f(u_{i} d_{i})= \left\{ \begin{array}{ll}
n+i-2,~~3\leq i \leq n-1, i~ \mbox{odd.}\\
n+i,~~2\leq i \leq n-2, i~ \mbox{even.}
\end{array} \right.
\end{displaymath}
Note that $f(u_1 d_1) = 2n-1$, and $f(u_n d_n) = n$. See the Figure~\ref{even-for-8} for the edge labeling.

Then the vertex sums(colors) $w(u_i)$ for $u_i, 1\leq i \leq n,$ are respectively $2n + \frac{n}{2} - 1$ for odd $i$, and $2n + \frac{n}{2} + 1$ for even $i$, except the case $i = n$ whose vertex sum is $n+1$. Notice that these two colors are distributed alternatively, and also distinct with those colors on pendent vertices since both are greater than $2n - 1$. On the other hand, the color $n+1$ of $u_n$ is repeated with that of $u_3$. Therefore the coloring we adopt here is a proper vertex coloring with $n+2$ colors, thus $\chi_{la} ( P_{n} \circ K_1 )\leq n+2$ and the proof is done.\qed
~~


\begin{lem}
$\chi_{la} ( P_{n} \circ K_1 )\leq n+2$ for $n \equiv 2~ (mod~ 4)$ and $n \ge 6$.
\end{lem}
\proof

Label the edges $u_{i} u_{i+1}, 1\leq i \leq n-1,$ as follows:

\begin{displaymath}
f(u_{i} u_{i+1})= \left\{ \begin{array}{ll}
\frac{n}{2} - \frac{i-1}{2},~~1\leq i \leq n-1, i~ \mbox{odd.}\\
n - \frac{i}{2},~~2\leq i \leq n-2, i~ \mbox{even.}
\end{array} \right.
\end{displaymath}

\begin{figure}[h]
\centering
\includegraphics[width=14cm]{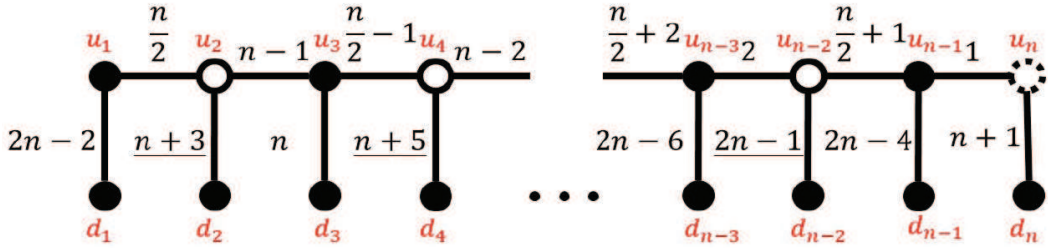}
\caption{The edge labeling for $n \equiv 2~(mod~4)$}
\label{even-for-10}
\end{figure}

Label the edges $u_{i} d_{i}, 1\leq i \leq n,$ as follows:

\begin{displaymath}
f(u_{i} d_{i})= \left\{ \begin{array}{ll}
n+i-3,~~3\leq i \leq n-1, i~ \mbox{odd.}\\
n+i+1,~~2\leq i \leq n-2, i~ \mbox{even.}
\end{array} \right.
\end{displaymath}
Note that $f(u_1 d_1) = 2n-2$, and $f(u_n d_n) = n+1$. See the Figure~\ref{even-for-10} for the edge labeling.

Then the vertex sums(colors) $w(u_i)$ for $u_i, 1\leq i \leq n,$ are respectively $2n + \frac{n}{2} - 2$ for odd $i$, and $2n + \frac{n}{2} + 2$ for even $i$, except for the case $i = n$ the vertex sum is $n+2$. Notice that these two colors are distributed alternatively, and also distinct with those colors on pendent vertices since both are greater than $2n - 1$. On the other hand the vertex sum(color) $n+2$ of $u_n$ is repeated with that of $u_5$. Therefore the coloring we adopt here is a proper vertex coloring with $n+2$ colors, thus $\chi_{la} ( P_{n} \circ K_1 )\leq n+2$ for $n \equiv 2~ (mod~ 4)$, $n \ge 6$, and the proof is done.\qed


\begin{lem}
$\chi_{la} ( P_{n} \circ K_1 )\leq n+2$ for $n \equiv 1~ (mod~ 2)$ and $n \ge 3$.
\end{lem}
\proof

Label the edges $u_{i} u_{i+1}, 1\leq i \leq n-1,$ as follows:

\begin{displaymath}
f(u_{i} u_{i+1})= \left\{ \begin{array}{ll}
n - \frac{i-1}{2},~~1\leq i \leq n-2, i~ \mbox{odd.}\\
\frac{n+1}{2} - \frac{i}{2},~~2\leq i \leq n-1, i~ \mbox{even.}
\end{array} \right.
\end{displaymath}

\begin{figure}[h]
\centering
\includegraphics[width=14cm]{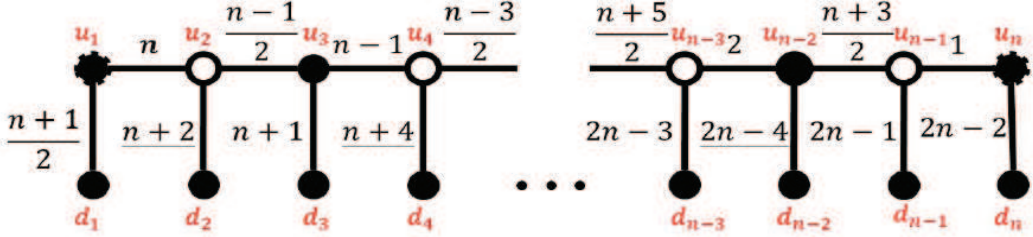}
\caption{The edge labeling for $n \equiv 1~(mod~2)$}
\label{odd}
\end{figure}

Label the edges $u_{i} d_{i}, 1\leq i \leq n,$ as follows:

\begin{displaymath}
f(u_{i} d_{i})= \left\{ \begin{array}{ll}
n+i-2,~~3\leq i \leq n, i~ \mbox{odd.}\\
n+i,~~2\leq i \leq n-1, i~ \mbox{even.}
\end{array} \right.
\end{displaymath}
Note that $f(u_1 d_1) = \frac{n+1}{2}$. See the Figure~\ref{odd} for the edge labeling for odd $n$.

Then the vertex sums(colors) $w(u_i)$ for $u_i, 1\leq i \leq n,$ are respectively $2n + \frac{n}{2} - \frac{1}{2}$ for odd $i$, except for the case $i = 1, n$, and $2n + \frac{n}{2} + \frac{3}{2}$ for even $i$. Note that first the vertex sum at $u_1$ is $n + \frac{n+1}{2}$, which is repeated with that of some $d_j$, if $n\ge 3$. Also note that the vertex sum at $u_n$ is $2n-1$, which is repeated with that of $d_{n-1}$. Therefore the coloring we adopt here is a proper vertex coloring with $n+2$ colors, thus $\chi_{la} ( P_{n} \circ K_1 )\leq n+2$ for odd $n$, $n \ge 3$, and the proof is done.\qed


\begin{lem}
$\chi_{la} ( P_{n} \circ K_1 )\geq n+2$ for $n \ge 4$.
\end{lem}
\proof

Note that the edges are labeled with $1,2,\cdots, 2n-1$.

Suppose $\chi_{la} ( P_{n} \circ K_1 )= n+1$ and there are $k$ of $n$ upper vertex sums(colors) $w(u_1), w(u_2), \cdots, w(u_n)$ are repeated with some of the lower vertex sums(colors) $w(d_1), w(d_2), \cdots, w(d_n)$. That is, there are $k$ of $u_i's$ for which $w(u_i) = w(d_j)$ for some $j$.

\begin{figure}[h]
\centering
\includegraphics[width=12cm]{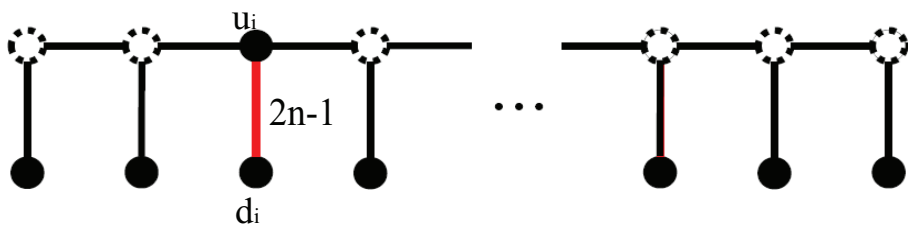}
\caption{There are $n-1$ repeated vertex sums(colors) for $P_{n} \circ K_1$}
\label{repeated-(n-1)}
\end{figure}

Note that the number of edges incident with the $k$ vertices with repeated upper vertex sums(colors) is $n+k-1$, since the number of edges not incident with these $k$ vertices is $n-k$. Let the total sum of the these $k$ repeated vertex sums be $\sigma$, then we see that $\sigma \leq k(2n-1)$ since each of the $k$ vertex sums are repeated with some $w(d_j)$ which is $\leq 2n-1$. On the other hand $\sigma \geq 1+2+\cdots +(n+k-1)$ if we use the smallest labels from those of edges incident with the $k$ vertices with repeated upper vertex sums(colors). Combining the above two inequalities we have $k(2n-1) \geq 1+2+\cdots +(n+k-1)$, which implies $k=n$ or $k=n-1$. However $k=n$ means $\chi_{la} ( P_{n} \circ K_1 )= n$, which is out of scope. Therefore it remains $k=n-1$, which we claim to reach a contradiction. The reason is that if $k=n-1$, we have only one non-repeated new vertex
sum(color) $w(u_i)$ at $u_i$ for some $i$, so that $f (u_i d_i ) = 2n - 1$ and
$w(d_i ) = 2n - 1$. See the Figure~\ref{repeated-(n-1)}.
\\
Note that there are two cases for $u_i$:\\

{\bf Case~1:} $deg(u_i)=2$. \\

In this case note that each of the $n-1$ repeated vertex sums(colors) is no greater $2n-1$. Look at the edge labeled $2n-2$, the only possibility is joining it with another edge labeled with $1$ and make another degree two vertex sum. Then next look at the edge labeled $2n-3$, the only possibility is joining it with another edge labeled with 1 or 2 respectively, otherwise it will be greater than $2n-1$. But one end vertex of the edge labeled $2n-3$ is of degree 3, a contradiction.\\

{\bf Case~2:} $deg(u_i)=3$. \\

Similarly in this situation we have two degree two vertex sums as $2n-2$ with 1 and $2n-3$ with 2. Then look at the edge labeled $2n-4$, the only way to place it is to join with two other edges labeled $1$ and $2$, but $1$ and $2$ are already used, a contradiction. Thus $k=n-1$ is impossible.

Note that $P_{n} \circ K_1$ is a tree with $n$ leaves. Thus $\chi_{la} ( P_{n} \circ K_1 )\geq n+2$ since previously it was proved $\chi_{la} (T) \geq l+1$ in case $T$ is a tree with $l$ leaves. Therefore we are done with the proof.\qed
~~

Therefore by previous results we have the following.
\begin{thm}
$\chi_{la} ( P_{n} \circ K_1 )= n+2$ for $n\ge 4$.
\end{thm}
%

\begin{figure}[h]
\centering
\includegraphics[width=12cm]{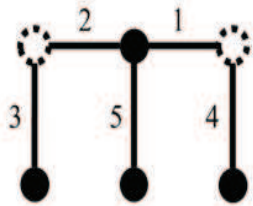}
\caption{The local antimagic chromatic number for $P_{3} \circ K_1$ is 4}
\label{P3-K1}
\end{figure}

Note that we have two exceptions $\chi_{la} ( P_{2} \circ K_1 )= 3$ and $\chi_{la} ( P_{3} \circ K_1 )= 4$(see Figure~\ref{P3-K1}).

That is, $P_{2} \circ K_1$ and $P_{3} \circ K_1$ belong to the class of trees $T$ such that $\chi_{la} ( T )= l+1$, where $l$ is the number of leaves of $T$, while $P_{n} \circ K_1$ belongs to the class of trees $T$ such that $\chi_{la} ( T )= l+2$ for $n\ge 4$.


\subsection{The Chromatic Number of $P_{n} \circ \overline{K_m}$}

Now we are in a position to proceed to the general case. Note that for $m \geq 2$ and $n \ge 2$ the graph $P_{n} \circ \overline{K_m}$ has $mn+n$ vertices, $mn+n-1$ edges, and $mn$ leaves. Generally we have the following.
\begin{thm}
$\chi_{la} ( P_{n} \circ \overline{K_m} )= mn+2$, for $m \geq 2, n \ge 2$.
\end{thm}

Let us start with the following lemma for the lower bound of the local antimagic chromatic number.

\begin{lem}\label{Pn-Km-lower-bound}
$\chi_{la} ( P_{n} \circ \overline{K_m} )\geq mn+2$ for $n \ge 2$ and $m\geq 2$.
\end{lem}
\proof

The graph $P_{n} \circ \overline{K_m}$ has $mn + n$ vertices and $mn + n -1$ edges. The edges are labeled with $1,2,\cdots, mn+n-1$. Note that $P_{n} \circ \overline{K_m}$ is a tree with $mn$ leaves, thus $\chi_{la} ( P_{n} \circ \overline{K_m} )\geq mn+1$ due to Theorem~\ref{leaf-lower-bound}. In order to show that $\chi_{la} ( P_{n} \circ \overline{K_m} )\geq mn+2$, it suffices to prove that $\chi_{la} ( P_{n} \circ \overline{K_m} )\neq mn+1$.

Suppose $\chi_{la} ( P_{n} \circ \overline{K_m} )= mn+1$ and there are exactly $k$ of $n$ upper vertex sums(colors) $w(u_1), w(u_2), \cdots, w(u_n)$ are repeated with some of the lower vertex sums(colors) $w(d_1), w(d_2), \cdots, w(d_{mn})$. That is, there are $k$ of $u_i's$ for which $w(u_i) = w(d_j)$ for some $j$. Note that in this situation one can only have an extra new color(vertex sum) other than the $mn$ lower vertex colors $w(d_1), w(d_2), \cdots, w(d_{mn})$ and no two adjacent vertices for $u_i$'s have the same color, thus it can be seen that $\lfloor \frac{n}{2} \rfloor \leq k \leq n-1$.

Note that the number of edges incident with the $k$ vertices with repeated upper vertex sums(colors) is $n+mk-1$, since the number of edges not incident with these $k$ vertices is $m(n-k)$. Let the total sum of the these $k$ repeated vertex sums be $\sigma$, then we see that $\sigma \leq k(n+mn-1)$ since each of the $k$ vertex sums are repeated with some $w(d_j)$ which is $\leq mn+n-1$.

On the other hand $\sigma \geq 1+2+\cdots +(n+mk-1)$ if we use the smallest labels from those of edges incident with the $k$ vertices with repeated upper vertex sums(colors). Combining the above two inequalities we have $k(n+mn-1) \geq 1+2+\cdots +(n+mk-1)$, which is a contradiction via the following calculation:
$$ 1+2+\cdots +(n+mk-1) - k(n+mn-1) $$
$$= \frac{1}{2}(n^2 + m^2 k^2 - n - mk - 2kn + 2k) $$
$$ = \frac{1}{2}[ ( n - \frac{1}{2} )^2 + ( mk - \frac{1}{2} )^2 - 2k(n-1) - \frac{1}{2}]$$
$$ > \frac{1}{2}[ ( n - 1 )^2 + ( mk - 1 )^2 - 2k(n-1) - \frac{1}{2}]$$
$$ = \frac{1}{2}[ ( n - 1 - k )^2 - k^2 + ( mk - 1 )^2 - \frac{1}{2}]$$
$$ = \frac{1}{2}[ ( n - 1 - k )^2 + (mk + k - 1)(mk - k -1) - \frac{1}{2}]$$
which is positive since $1 \leq \lfloor \frac{n}{2} \rfloor \leq k \leq n-1$ and $m \geq 2$. Therefore we are done with the proof.\qed
~~

To obtain the upper bound for the chromatic number $\chi_{la} ( P_{n} \circ \overline{K_m} )$, we give the the following first for $m=2$:

\begin{thm}
$\chi_{la} ( P_{n} \circ \overline{K_2} )= 2n+2$, for $n \ge 2$.
\end{thm}

We need the following lemmas in oder to obtain the above result.
\begin{lem}\label{Pn-K2-even}
$\chi_{la} ( P_{n} \circ \overline{K_2} )\leq 2n+2$ for $n \equiv 0~ (mod~ 2)$ and $n \ge 2$.
\end{lem}

\proof

Assume the vertex set is $V(P_{n} \circ \overline{K_2}) = \{ m_1,m_2, \cdots, m_n, u_1,u_2, \cdots, u_n, d_1,d_2, \cdots, d_n \}$ and the edge set is $E(P_{n} \circ \overline{K_2}) = \{ m_i m_{i+1} : ~ 1 \leq i \leq n-1 \} \cup \{ m_i u_i : ~1 \leq i \leq n \}\cup \{ m_i d_i : ~1 \leq i \leq n \}$. See Figure~\ref{Pn-K2-even}.

\begin{figure}[h]
\centering
\includegraphics[width=14cm]{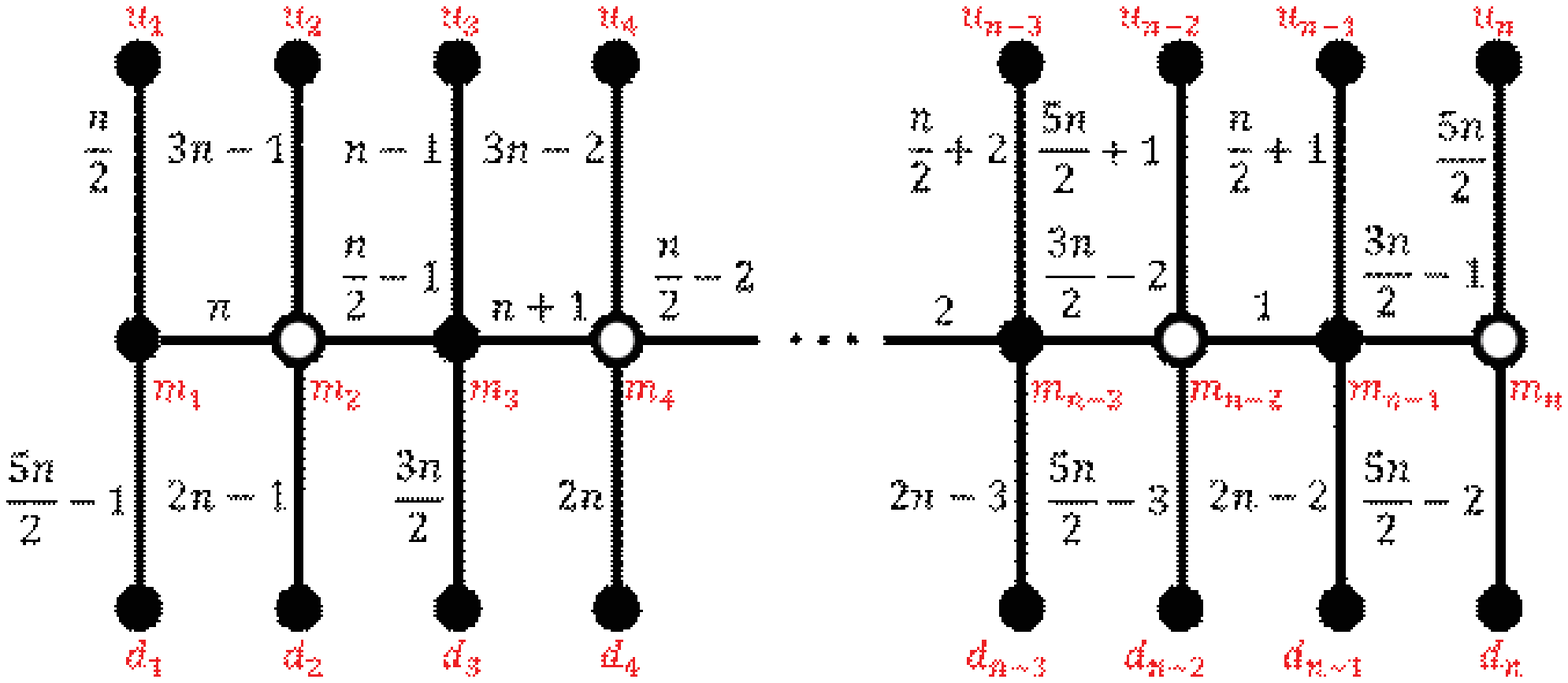}
\caption{The edge labeling for $P_{n} \circ \overline{K_2}, n \equiv 0~(mod~2)$}
\label{Pn-K2-even}
\end{figure}

Label the edges $m_{i} m_{i+1}, 1\leq i \leq n-1,$ as follows:

\begin{displaymath}
f(m_{i} m_{i+1})= \left\{ \begin{array}{ll}
n + \frac{i-1}{2},~~1\leq i \leq n-1, i~ \mbox{odd.}\\
\frac{n}{2} - \frac{i}{2},~~2\leq i \leq n-2, i~ \mbox{even.}
\end{array} \right.
\end{displaymath}

Label the edges $m_{i} u_{i}, 1\leq i \leq n,$ as follows:

\begin{displaymath}
f(m_{i} u_{i})= \left\{ \begin{array}{ll}
\frac{n}{2}, ~~i=1.\\
n - \frac{i-1}{2},~~3\leq i \leq n-1, i~ \mbox{odd.}\\
3n - \frac{i}{2},~~2\leq i \leq n, i~ \mbox{even.}
\end{array} \right.
\end{displaymath}

Label the edges $m_{i} d_{i}, 1\leq i \leq n,$ as follows:

\begin{displaymath}
f(m_{i} d_{i})= \left\{ \begin{array}{ll}
\frac{5n}{2} - 1, ~~i=1.\\
\frac{3n}{2} - \frac{i-3}{2},~~3\leq i \leq n-1, i~ \mbox{odd.}\\
2n + \frac{i-4}{2},~~2\leq i \leq n, i~ \mbox{even.}
\end{array} \right.
\end{displaymath}

Then the vertex sums(colors) $w(m_i)$ for $m_i, 1\leq i \leq n,$ are respectively $4n-1$ for odd $i$, and $\frac{13n}{2} -3 $ for even $i$. Therefore the coloring we adopt here is a proper vertex coloring with $2n+2$ colors, thus the proof is done.\qed
~~

\begin{lem}\label{Pn-K2-odd}
$\chi_{la} ( P_{n} \circ \overline{K_2} )\leq 2n+2$ for $n \equiv 1~ (mod~ 2)$ and $n \ge 3$.
\end{lem}

\proof

Again assume the vertex set is $V(P_{n} \circ \overline{K_2}) = \{ m_1,m_2, \cdots, m_n, u_1,u_2, \cdots, u_n, d_1,d_2, \cdots, d_n \}$ and the edge set is $E(P_{n} \circ \overline{K_2}) = \{ m_i m_{i+1} : ~ 1 \leq i \leq n-1 \} \cup \{ m_i u_i : ~1 \leq i \leq n \}\cup \{ m_i d_i : ~1 \leq i \leq n \}$. See Figure~\ref{Pn-K2-odd}.

\begin{figure}[h]
\centering
\includegraphics[width=14cm]{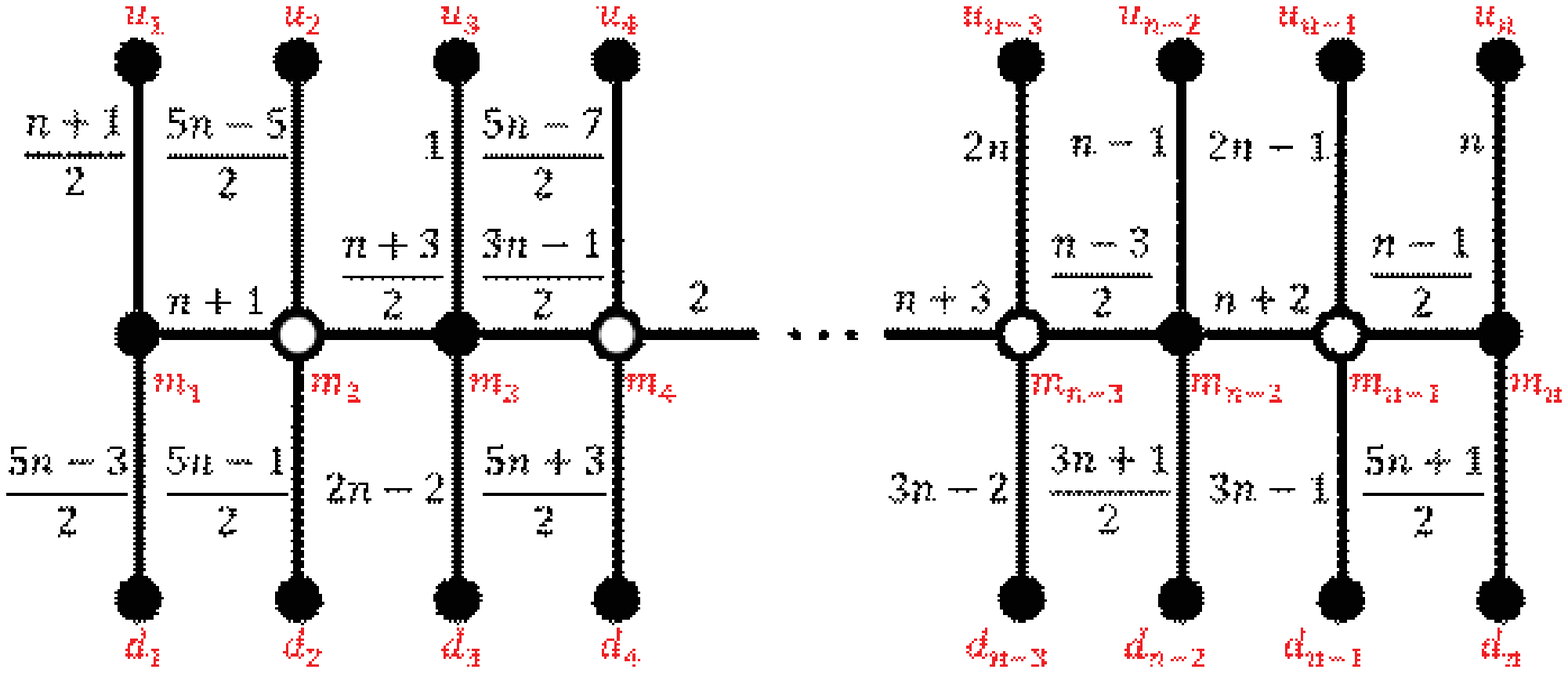}
\caption{The edge labeling for $P_{n} \circ \overline{K_2}, n \equiv 1~(mod~2)$}
\label{Pn-K2-odd}
\end{figure}

Label the edges $m_{i} m_{i+1}, 1\leq i \leq n-1,$ as follows:

\begin{displaymath}
f(m_{i} m_{i+1})= \left\{ \begin{array}{ll}
n + 1,~~i=1.\\
\frac{n+3}{2},~~i=2.\\
\frac{3n}{2} - \frac{i-2}{2},~~3\leq i \leq n, i~ \mbox{odd.}\\
\frac{i}{2},~~4\leq i \leq n-1, i~ \mbox{even.}
\end{array} \right.
\end{displaymath}

Label the edges $m_{i} u_{i}, 1\leq i \leq n,$ as follows:

\begin{displaymath}
f(m_{i} u_{i})= \left\{ \begin{array}{ll}
\frac{n+1}{2},~~i=1.\\
1,~~i=3.\\
\frac{n}{2} + \frac{i}{2},~~5\leq i \leq n, i~ \mbox{odd.}\\
\frac{5n}{2} - \frac{i+3}{2},~~2\leq i \leq n, i~ \mbox{even.}
\end{array} \right.
\end{displaymath}

Label the edges $m_{i} d_{i}, 1\leq i \leq n,$ as follows:

\begin{displaymath}
f(m_{i} d_{i})= \left\{ \begin{array}{ll}
\frac{5n-3}{2}, ~~i=1.\\
\frac{5n-1}{2},~~i=2.\\
2n - \frac{i+1}{2},~~3\leq i \leq n-2, i~ \mbox{odd.}\\
\frac{5n}{2} + \frac{i-1}{2},~~4\leq i \leq n-1, i~ \mbox{even.}\\
\frac{5n+1}{2},~~i=n.\\
\end{array} \right.
\end{displaymath}

Then the vertex sums(colors) $w(m_i)$ for $m_i, 1\leq i \leq n,$ are respectively $4n$ for odd $i$, and $\frac{13n-1}{2}$ for even $i$. Therefore the coloring we adopt here is a proper vertex coloring with $2n+2$ colors, thus the proof is done.\qed
~~

Next we proceed to calculate the general case $P_{n} \circ \overline{K_m}$ for $m \geq 3$. Let us start with $m$ even, $m \geq 4$ as in the following.

\begin{lem}
$\chi_{la} ( P_{n} \circ \overline{K_m} )\leq mn+2$ for $m$ even, $m \geq 4$ and $n \ge 2$.
\end{lem}
\proof
Assume the vertex set is $V(P_{n} \circ \overline{K_m}) = \bigcup_{i=1}^{n}\{ d_i^1,d_i^2, \cdots, d_i^m\}\cup\{ u_1,u_2, \cdots, u_n \}$ and the edge set is $E(P_{n} \circ \overline{K_m}) = \bigcup_{i=1}^{n} \{ u_i d_i^j : ~1 \leq i \leq n, 1 \leq j \leq m \}\cup\{ u_i u_{i+1} : ~ 1 \leq i \leq n-1 \}$.

With labels used in Lemma~\ref{Pn-K2-even} and Lemma~\ref{Pn-K2-odd}, we give the labeling $1,2,\cdots, 3n-1$ to the edges $u_i u_{i+1}, 1 \leq i \leq n-1$ and $u_i d_i^j : ~1 \leq i \leq n, 1 \leq j \leq 2$ for the subgraph $P_{n} \circ \overline{K_2}$ first, so that the partial vertex sums at $u_i$ is $4n-1$ for odd $i$, and $\frac{13n}{2} -3 $ for even $i$ whenever $n$ is even, and $4n$ for odd $i$, and $\frac{13n-1}{2}$ for even $i$ whenever $n$ is odd, respectively.

Then we give the edge labels $3n, 3n+1, \cdots, mn+n-1$ to the remaining edges as follows whenever $1\leq i \leq n$ and $3\leq j \leq m$:

\begin{displaymath}
f(u_{i} d_i^{j})= \left\{ \begin{array}{ll}
jn -1 + i, ~~j~ \mbox{odd}.\\
(j+1)n - i,~~j~ \mbox{even}.
\end{array} \right.
\end{displaymath}
Therefore it can be calculated the vertex sums at $u_i$ is $(\frac{m^2 + 2m}{2})n - \frac{m}{2}$ for odd $i$, and $(\frac{m^2 + 2m +5}{2})n - \frac{m+4}{2}$ for even $i$ whenever $n$ is even, and $(\frac{m^2 + 2m}{2})n - \frac{m-2}{2}$ for odd $i$, and $(\frac{m^2 + 2m +5}{2})n - \frac{m-1}{2}$ for even $i$ whenever $n$ is odd, respectively. Thus we obtain two new colors(vertex sums) other than those over leaves. We are done with the proof.\qed
~~

Then we consider the case $m$ odd, $m \geq 3$ as in the following. There are two subcases for $n$ odd and $n$ even respectively.

\begin{lem}
$\chi_{la} ( P_{n} \circ \overline{K_m} )\leq mn+2$ for $m, n$ odd, $m \geq 3$ and $n \ge 3$.
\end{lem}
\proof
Assume the vertex set is $V(P_{n} \circ \overline{K_m}) = \bigcup_{i=1}^{n}\{ d_i^1,d_i^2, \cdots, d_i^m\}\cup\{ u_1,u_2, \cdots, u_n \}$ and the edge set is $E(P_{n} \circ \overline{K_m}) = \bigcup_{i=1}^{n} \{ u_i d_i^j : ~1 \leq i \leq n, 1 \leq j \leq m \}\cup\{ u_i u_{i+1} : ~ 1 \leq i \leq n-1 \}$.

We first give the labeling $1,2,\cdots, 4n-1$ to the edges $u_i u_{i+1}, 1 \leq i \leq n-1$ and $u_i d_i^j : ~1 \leq i \leq n, 1 \leq j \leq 3$ for the subgraph $P_{n} \circ \overline{K_3}$ as in the following.

Label the edges $u_{i} u_{i+1}, 1 \leq i \leq n-1$:
\begin{displaymath}
f(u_{i} u_{i+1})= \left\{ \begin{array}{ll}
n-1 - \frac{i-1}{2}, ~~i~ \mbox{odd}.\\
\frac{i}{2},~~i~ \mbox{even}.
\end{array} \right.
\end{displaymath}

Label the edges $u_i d_i^j : ~1 \leq i \leq n, 1 \leq j \leq 3$:
\begin{displaymath}
f(u_i d_i^1)= \left\{ \begin{array}{ll}
\frac{3n-1}{2} + \frac{i-1}{2}, ~~i~ \mbox{odd}, i\neq n.\\
n + \frac{i-1}{2},~~i~ \mbox{even}.\\
\frac{5n-3}{2},~~i=n.
\end{array} \right.
\end{displaymath}
\begin{displaymath}
f(u_i d_i^2)= \left\{ \begin{array}{ll}
\frac{5n-1}{2} + \frac{i-1}{2}, ~~i~ \mbox{odd}.\\
2n-1 + \frac{i-2}{2},~~i~ \mbox{even}.
\end{array} \right.
\end{displaymath}
\begin{displaymath}
f(u_i d_i^3)= \left\{ \begin{array}{ll}
4n-i, ~~i~ \mbox{odd}.\\
4n-i,~~i~ \mbox{even}.
\end{array} \right.
\end{displaymath}
So the partial vertex sums at $u_i$ is $9n-3$ for odd $i$, and $8n -3 $ for even $i$ respectively.

Then secondly we give the edge labels $4n, 4n+1, \cdots, mn+n-1$ to the remaining edges as follows whenever $1\leq i \leq n$ and $4\leq j \leq m$:

\begin{displaymath}
f(u_{i} d_i^{j})= \left\{ \begin{array}{ll}
(j+1)n - i, ~~j~ \mbox{odd}.\\
jn - 1 + i,~~j~ \mbox{even}.
\end{array} \right.
\end{displaymath}
Therefore it can be calculated the vertex sums at $u_i$ is $(\frac{m^2 + 2m + 3}{2})n - \frac{m+3}{2}$ for odd $i$, and $(\frac{m^2 + 2m +1}{2})n - \frac{m+3}{2}$ for even $i$, respectively. Thus we obtain two new colors(vertex sums) other than those over leaves. We are done with the proof.\qed
~~

\begin{lem}
$\chi_{la} ( P_{n} \circ \overline{K_m} )\leq mn+2$ for $m$ odd and $n$ even, $m \geq 3$ and $n \ge 4$.
\end{lem}
\proof
Assume the vertex set is $V(P_{n} \circ \overline{K_m}) = \bigcup_{i=1}^{n}\{ d_i^1,d_i^2, \cdots, d_i^m\}\cup\{ u_1,u_2, \cdots, u_n \}$ and the edge set is $E(P_{n} \circ \overline{K_m}) = \bigcup_{i=1}^{n} \{ u_i d_i^j : ~1 \leq i \leq n, 1 \leq j \leq m \}\cup\{ u_i u_{i+1} : ~ 1 \leq i \leq n-1 \}$.

We first give the labeling $1,2,\cdots, 4n-1$ to the edges $u_i u_{i+1}, 1 \leq i \leq n-1$ and $u_i d_i^j : ~1 \leq i \leq n, 1 \leq j \leq 3$ for the subgraph $P_{n} \circ \overline{K_3}$ as in the following.

Label the edges $u_{i} u_{i+1}, 1 \leq i \leq n-1$:
\begin{displaymath}
f(u_{i} u_{i+1})= \left\{ \begin{array}{ll}
n-1 - \frac{i-1}{2}, ~~i~ \mbox{odd}.\\
\frac{i}{2},~~i~ \mbox{even}.
\end{array} \right.
\end{displaymath}

Label the edges $u_i d_i^j : ~1 \leq i \leq n, 1 \leq j \leq 3$:
\begin{displaymath}
f(u_i d_i^1)= \left\{ \begin{array}{ll}
n + \frac{i-1}{2}, ~~i~ \mbox{odd}.\\
2n + \frac{i-2}{2},~~i~ \mbox{even}, i\neq n.\\
\frac{5n}{2},~~i=n.
\end{array} \right.
\end{displaymath}
\begin{displaymath}
f(u_i d_i^2)= \left\{ \begin{array}{ll}
\frac{3n}{2} + \frac{i-1}{2}, ~~i~ \mbox{odd}.\\
\frac{7n}{2}-i,~~i~ \mbox{even}, i\neq n.\\
\frac{7n}{2}-1, ~~i=n.
\end{array} \right.
\end{displaymath}
\begin{displaymath}
f(u_i d_i^3)= \left\{ \begin{array}{ll}
\frac{7n}{2} -2 -i, ~~i~ \mbox{odd}.\\
\frac{7n}{2} + \frac{i}{2},~~i~ \mbox{even}, i\neq n.\\
\frac{7n}{2}, ~~i=n.
\end{array} \right.
\end{displaymath}
So the partial vertex sums at $u_i$ is $7n-4$ for odd $i$, and $10n -1$ for even $i$ respectively.

Then secondly we give the edge labels $4n, 4n+1, \cdots, mn+n-1$ to the remaining edges as follows whenever $1\leq i \leq n$ and $4\leq j \leq m$:

\begin{displaymath}
f(u_{i} d_i^{j})= \left\{ \begin{array}{ll}
(j+1)n - i, ~~j~ \mbox{odd}.\\
jn - 1 + i,~~j~ \mbox{even}.
\end{array} \right.
\end{displaymath}
Therefore it can be calculated the vertex sums at $u_i$ is $(\frac{m^2 + 2m - 1}{2})n - \frac{m+5}{2}$ for odd $i$, and $(\frac{m^2 + 2m +5}{2})n - \frac{m-1}{2}$ for even $i$, respectively. Thus we obtain two new colors(vertex sums) other than those over leaves. We are done with the proof.\qed
~~
Therefore we conclude with the following.
\begin{thm}
$\chi_{la} ( P_{n} \circ \overline{K_m} )= mn+2$, for $m \geq 2, n \ge 2$.
\end{thm}

\section{Corona Product of Cycles with Null Graphs}

\subsection{The Chromatic Number of $C_{n} \circ K_1$}

\begin{figure}[h]
\centering
\includegraphics[width=14cm]{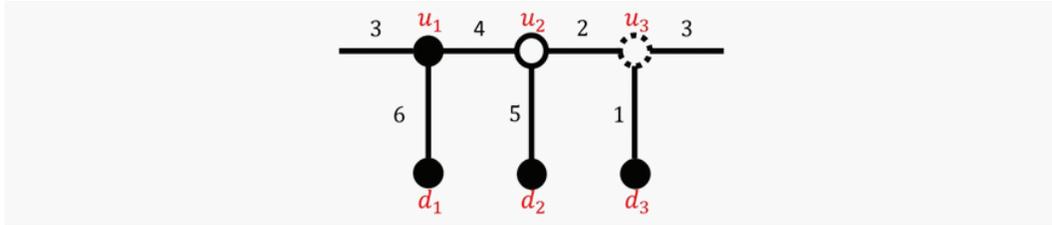}
\caption{$\chi_{la} ( C_{3} \circ K_1 )= 5$}
\label{C3-K1}
\end{figure}

It is not hard to calculate the coloring number for $C_{3} \circ K_1 $. See Figure~\ref{C3-K1}
\begin{lem}\label{Cn_K1_odd_n}
$\chi_{la} ( C_{3} \circ K_1 )= 5$.
\end{lem}

Then we proceed with general cases.

\begin{lem}\label{Cn_K1_odd_n}
$\chi_{la} ( C_{n} \circ K_1 )\leq n+2$ for $n \equiv 1~ (mod~ 2)$ and $n \ge 5$.
\end{lem}
\proof

Label the edges $u_{i} u_{i+1}, 1\leq i \leq n-1,$ as follows:

\begin{displaymath}
f(u_{i} u_{i+1})= \left\{ \begin{array}{ll}
n - \frac{i-1}{2},~~1\leq i \leq n-2, i~ \mbox{odd.}\\
\frac{n+1}{2} - \frac{i}{2},~~2\leq i \leq n-1, i~ \mbox{even.}
\end{array} \right.
\end{displaymath}

\begin{figure}[h]
\centering
\includegraphics[width=14cm]{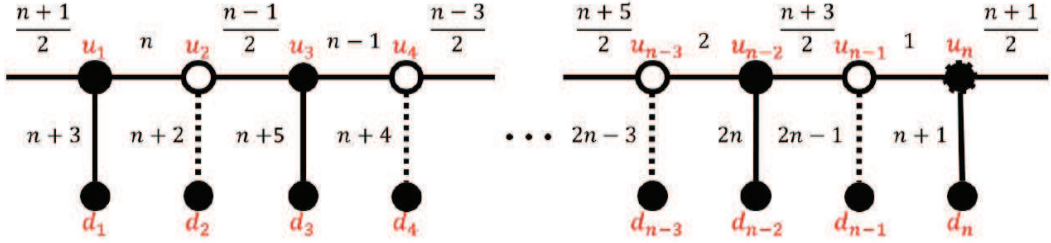}
\caption{The edge labeling for $n \equiv 1~(mod~2)$}
\label{Cn-odd}
\end{figure}

Note that $f(u_n u_1) = \frac{n+1}{2}$. Then label the edges $u_{i} d_{i}, 1\leq i \leq n,$ as follows:

\begin{displaymath}
f(u_{i} d_{i})= \left\{ \begin{array}{ll}
n+i+2,~~1\leq i \leq n-2, i~ \mbox{odd.}\\
n+i,~~2\leq i \leq n-1, i~ \mbox{even.}
\end{array} \right.
\end{displaymath}
Note that $f(u_n d_n) = n+1$. See the Figure~\ref{Cn-odd} for the edge labeling for odd $n$.

Then the vertex sums(colors) $w(u_i)$ for $u_i, 1\leq i \leq n,$ are respectively $\frac{5n+7}{2}$ for odd $i$, except for the case $i = n$, and $\frac{5n+3}{2}$ for even $i$. Note that the vertex sum at $u_n$ is $\frac{3n+5}{2}$, which is not repeated with that of $u_{n-1}$ and $d_n$ for $n \ge 5$. Also notice that the vertex sum at $u_n$ is the same with that of some $d_j$. Therefore the coloring we adopt here is a proper vertex coloring with $n+2$ colors, thus $\chi_{la} ( C_{n} \circ K_1 )\leq n+2$ for $n \equiv 1~ (mod~ 2)$, $n \ge 5$, and the proof is done.\qed
~~


\begin{lem}
$\chi_{la} ( C_{n} \circ K_1 )\leq n+2$ for $n \equiv 0~ (mod~ 2)$ and $n \ge 6$.
\end{lem}
\proof

Label the edges $u_{i} u_{i+1}, 1\leq i \leq n-1,$ as follows:

\begin{displaymath}
f(u_{i} u_{i+1})= \left\{ \begin{array}{ll}
\frac{n + i-3}{2},~~1\leq i \leq n-1, i~ \mbox{odd.}\\
n - \frac{i-2}{2},~~2\leq i \leq n-2, i~ \mbox{even.}
\end{array} \right.
\end{displaymath}

\begin{figure}[h]
\centering
\includegraphics[width=14cm]{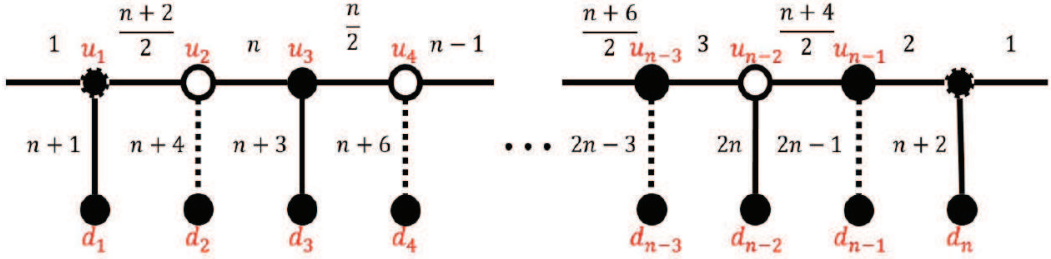}
\caption{The edge labeling for $n \equiv 0~ (mod~2)$}
\label{Cn-even}
\end{figure}

Note that $f(u_n u_1) = 1$. Then label the edges $u_{i} d_{i}, 1\leq i \leq n,$ as follows:

\begin{displaymath}
f(u_{i} d_{i})= \left\{ \begin{array}{ll}
n+i,~~1\leq i \leq n-2, i~ \mbox{odd.}\\
n+i+2,~~2\leq i \leq n-1, i~ \mbox{even.}
\end{array} \right.
\end{displaymath}
Note that $f(u_n d_n) = n+2$. See the Figure~\ref{Cn-even} for the edge labeling for odd $n$.

Then the vertex sums(colors) $w(u_i)$ for $u_i, 1\leq i \leq n,$ are respectively $\frac{5n+6}{2}$ for odd $i$, except for the case $i = 1$, and $\frac{5n+10}{2}$ for even $i$. Note that the vertex sum at $u_1$ is $\frac{3n+6}{2}$, which is not repeated with that of $u_{2}$ and $d_1$ for $n \ge 6$. Also notice that the vertex sum at $u_1$ is the same with that of some $d_j$ since $n \ge 6$. On the other hand, $w(u_n) = n+5$ is the same with $w(d_5)$. Therefore the coloring we adopt here is a proper vertex coloring with $n+2$ colors, thus $\chi_{la} ( C_{n} \circ K_1 )\leq n+2$ for even $n$, $n \ge 6$, and the proof is done.\qed
~~


\begin{lem}
$\chi_{la} ( C_{n} \circ K_1 )\geq n+2$ for $n \ge 4$.
\end{lem}
\proof

Note that there are $2n$ edges for the graph $C_{n} \circ K_1 $, and the local antimagic chromatic number $\chi_{la} ( C_{n} \circ K_1 )\geq n$ since it has $n$ pendent vertices with different vertex sums(colors).

First we suppose that $\chi_{la} ( C_{n} \circ K_1 )= n$. Then may see that the vertex sums(colors) $w(u_1), w(u_2), \cdots, w(u_n)$ are repeated with some $w(d_j)$'s so that $w(u_i) \leq 2n$ for each $i$, which is impossible. The reason is that, for the particular vertex $u_i$ incident with the edge labeled $2n$, the vertex sum $w(u_i)$ must be greater than $2n$.

Then consider the case $\chi_{la} ( C_{n} \circ K_1 )= n+1$ and suppose that there are $k$ of $n$ upper vertex sums(colors) $w(u_1), w(u_2), \cdots, w(u_n)$ are repeated with some of the lower vertex sums(colors) $w(d_1), w(d_2), \cdots, w(d_n)$. That is, there are $k$ of $u_i$'s for which $w(u_i) = w(d_j)$ for some $j$. Note that the number of edges incident with the $k$ vertices with repeated upper vertex sums(colors) is $n+k$, since the number of edges not incident with these $k$ vertices is $n-k$. Let the total sum of the these $k$ vertex sums be $\sigma$, then we see that $\sigma \leq k(2n)$ since each of the $k$ vertex sums are repeated with some $d_j$ which is $\leq 2n$. On the other hand $\sigma \geq 1+2+\cdots +(n+k)$ if we use the smallest labels from those of edges incident with the $k$ vertices with repeated upper vertex sums(colors). Combining the above two inequalities we have $k(2n) \geq 1+2+\cdots +(n+k)$, which implies $(n-k)^2 + n + k \leq 0$, a contradiction.

Therefore $\chi_{la} ( C_{n} \circ K_1 )\geq n+2$ for $n \ge 4$, the proof is done.\qed
~~

\begin{figure}[h]
\centering
\includegraphics[width=14cm]{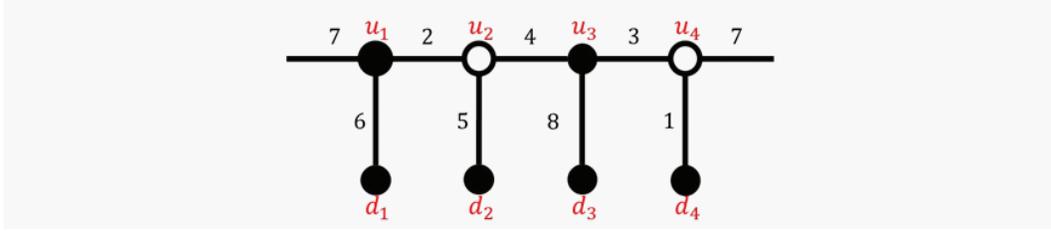}
\caption{The local antimagic coloring number $\chi_{la}(C_{4} \circ {K_1}) \leq 6$}
\label{C4-K1}
\end{figure}

\begin{thm}
$\chi_{la} ( C_{n} \circ K_1 )= n+2$ for $n\ge 4$.
\end{thm}
\proof
By previous lemmas and the Figure~\ref{C4-K1}. \qed
~~

\subsection{The Chromatic Number of $C_{n} \circ \overline{K_m}$}

Now we are in a position to proceed to the general case. Note that for $m \geq 2$ and $n \ge 3$ the graph $C_{n} \circ \overline{K_m}$ has $mn+n$ vertices, $mn+n$ edges, and $mn$ leaves.

With the similar method in Lemma~\ref{Pn-Km-lower-bound}, we have the following lemma for the lower bound of the local antimagic chromatic number.

\begin{figure}[h]
\centering
\includegraphics[width=14cm]{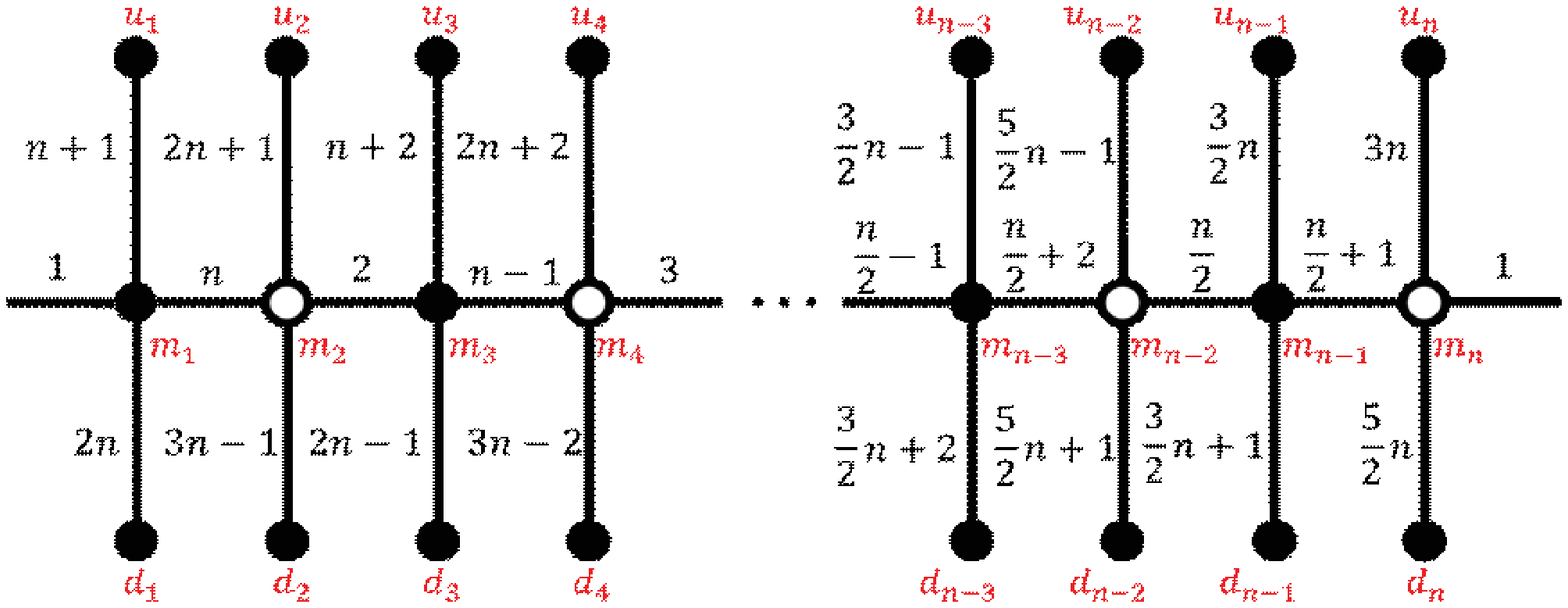}
\caption{The edge labeling for $C_{n} \circ \overline{K_2}$, $n \equiv 0~ (mod~2)$}
\label{Cn-K2-even}
\end{figure}

\begin{lem}\label{Cn-Km-lower-bound}
$\chi_{la} ( C_{n} \circ \overline{K_m} )\geq mn+2$ for $n \ge 3$ and $m\geq 2$.
\end{lem}
Then we turn to the upper bounds of the chromatic numbers. Let us start with the graph $C_{n} \circ \overline{K_2}$(see Figure~\ref{Cn-K2-even}). Suppose the vertex set is $V(C_{n} \circ \overline{K_2}) = \{ m_1, \cdots, m_n \} \cup \{ u_1, \cdots, u_n \} \cup\{ d_1, \cdots, d_n \}$ where $\{ m_1, \cdots, m_n \}$ are vertices of degree four inducing the cycle $C_n$, and $u_i, d_i$ are pendent vertices for each $1\leq i \leq n$. The edge set is $\{ m_i m_{i+1}:~ 1\leq i \leq n-1\} \cup \{ m_n m_1 \} \cup \{ m_i u_i: ~1\leq i \leq n \} \cup \{ m_i d_i: ~1\leq i \leq n \}$.

\begin{lem}\label{Cn-K2-n-even}
$\chi_{la} ( C_{n} \circ \overline{K_2} )\leq 2n+2$ for $n \equiv 0~ (mod~ 2)$ and $n \ge 4$.
\end{lem}
\proof

Label the edges over the cycle $m_{i} m_{i+1}, 1\leq i \leq n-1,$ as follows:

\begin{displaymath}
f(m_{i} m_{i+1})= \left\{ \begin{array}{ll}
n - \frac{i -1}{2}, ~i~ \mbox{odd.}\\
2 + \frac{i-2}{2}, ~i~ \mbox{even.}
\end{array} \right.
\end{displaymath}

Also label $f(m_n m_1) = 1$.

Then label the edges $m_{i} u_{i}, 1\leq i \leq n,$ as follows:

\begin{displaymath}
f(m_i u_i)= \left\{ \begin{array}{ll}
n+1+\frac{i-1}{2},~~ i~ \mbox{odd.}\\
2n+1+\frac{i-2}{2},~~i~ \mbox{even,} i\ne n.\\
3n, ~~ i=n.
\end{array} \right.
\end{displaymath}

Label the edges $m_{i} d_{i}, 1\leq i \leq n,$ as follows:

\begin{displaymath}
f(m_i d_i)= \left\{ \begin{array}{ll}
2n-\frac{i-1}{2},~~ i~ \mbox{odd.}\\
3n-1-\frac{i-2}{2},~~i~ \mbox{even.}
\end{array} \right.
\end{displaymath}

Then the vertex sums(colors) $w(m_i)$ for $m_i, 1\leq i \leq n,$ are respectively $4n+2$ for odd $i$, and $6n+2$ for even $i$.
Therefore the coloring we adopt here is a proper vertex coloring with $2n+2$ colors, thus $\chi_{la} ( C_{n} \circ \overline{K_2} )\leq 2n+2$ and we are done.\qed
~~

\begin{lem}
$\chi_{la} ( C_{n} \circ \overline{K_2} )\leq 2n+2$ for $n \equiv 1~ (mod~ 2)$ and $n \ge 5$.
\end{lem}
\proof

Label the edges over the cycle $m_{i} m_{i+1}, 1\leq i \leq n-1,$ as follows:

\begin{displaymath}
f(m_{i} m_{i+1})= \left\{ \begin{array}{ll}
n+1 - \frac{i -1}{2}, i~ \mbox{odd.}\\
3 + \frac{i-2}{2}, i~ \mbox{even.}
\end{array} \right.
\end{displaymath}

Also label $f(m_n m_1) = 1$.

\begin{figure}[h]
\centering
\includegraphics[width=14cm]{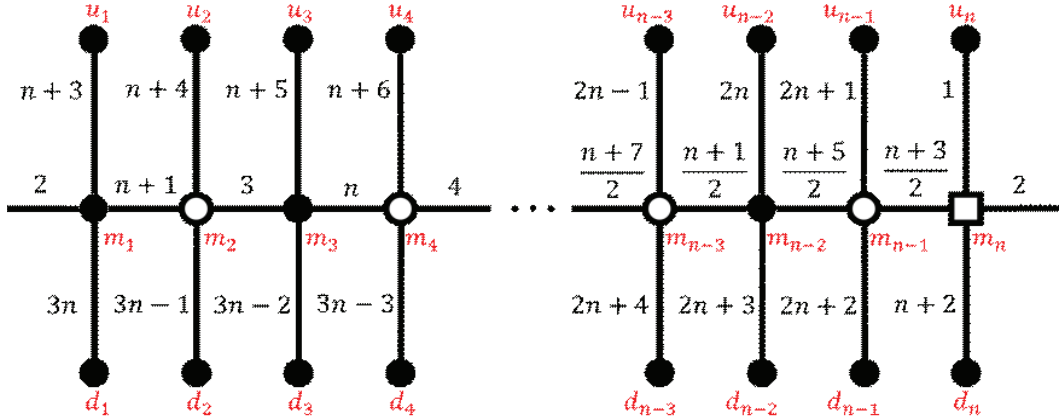}
\caption{The edge labeling for $C_{n} \circ \overline{K_2}$, $n \equiv 1~ (mod~2)$}
\label{Cn-K2-odd}
\end{figure}

Then label the edges $m_{i} u_{i}, 1\leq i \leq n,$ as follows:

\begin{displaymath}
f(m_i u_i)= \left\{ \begin{array}{ll}
n+2+i,~~ 1\leq i \leq n-1.\\
1, ~~i=n.
\end{array} \right.
\end{displaymath}

Label the edges $m_{i} d_{i}, 1\leq i \leq n,$ as follows:

\begin{displaymath}
f(m_i d_i)= \left\{ \begin{array}{ll}
3n+1-i,~~ 1\leq i \leq n-1.\\
n+2, ~~i=n.
\end{array} \right.
\end{displaymath}

Then the vertex sums(colors) $w(m_i)$ for $m_i, 1\leq i \leq n,$ are respectively $5n+6$ for odd $i$ except $i=n$, and $5n+7$ for even $i$. Note that $w(m_n)=\frac{3n+13}{2}$, which is between $n+3$ and $3n$ whenever $n \ge 5$ and is repeated.
Therefore it can be easily checked that the vertex coloring we adopt here is a proper one with $2n+2$ colors, thus $\chi_{la} ( C_{n} \circ \overline{K_2} )\leq 2n+2$ and we are done.\qed
~~

Combining previous lemmas we may obtain the following fact for $C_{n} \circ \overline{K_2}$.

\begin{thm}
$\chi_{la} ( C_{n} \circ \overline{K_2} )= 2n+2$ for $n \ge 4$.
\end{thm}

Note that we have the following exception for $\chi_{la} ( C_{n} \circ \overline{K_2} )$ as $n=3$.
\begin{lem}
$\chi_{la} ( C_{3} \circ \overline{K_2} )= 9$.
\end{lem}
\begin{figure}[h]
\centering
\includegraphics[width=7cm]{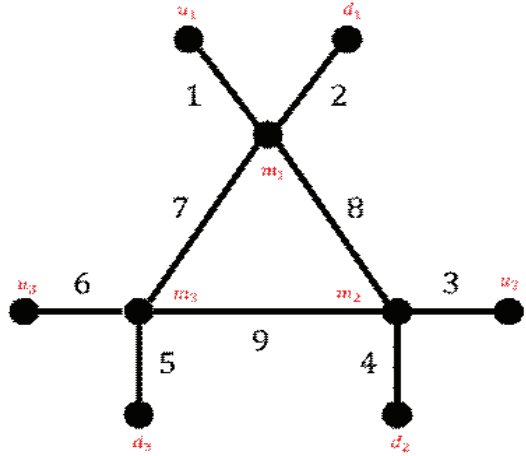}
\caption{The coloring number of $C_3 \circ K_2$ is 9}
\label{C3-K2}
\end{figure}
\proof
Assume the vertices of $C_3 \circ K_2$ are $m_1,m_2,m_3,u_1,u_2,u_3,d_1,d_2,d_3$ as in the Figure~\ref{C3-K2}. Note that the colors at $m_1,m_2,m_3$ must be distinct. Suppose $\chi_{la} ( C_{3} \circ \overline{K_2} )$ is 8 or less. Hence the color at some $m_i$ is repeated with those over leaves $u_j$ and $d_k$. We see that the color(vertex sum) at $m_i \geq 1+2+3+4 = 10$ with the smallest possible edge labels, however the largest edge label is 9, a contradiction. Therefore the chromatic number $\chi_{la} ( C_{3} \circ \overline{K_2} )$ is at least 9. \qed
~~

In fact generally we have the following for $C_{3} \circ \overline{K_m}, m \geq 2$.

\begin{thm}\label{C3-Km}
$\chi_{la} ( C_{3} \circ \overline{K_m} )= 3m+3$, for $m\geq 2$.
\end{thm}
\proof
Again assume the vertices of the subgraph $C_3$ of $C_3 \circ K_2$ are $m_1,m_2,m_3$. Note that the colors at $m_1,m_2,m_3$ must be distinct. Suppose $\chi_{la} ( C_{3} \circ \overline{K_m} )$ is $3m+2$ or less. Hence the color at some $m_i$ must be repeated with those over leaves. We see that the color(vertex sum) at $m_i \geq 1+2+\cdots+(m+2) = \frac{(m+3)(m+2)}{2}$ with the smallest possible edge labels, however the largest edge label is $3m+3$, a contradiction since $\frac{(m+3)(m+2)}{2} > 3m+3$ for $m \geq 2$. Therefore the chromatic number $\chi_{la} ( C_{3} \circ \overline{K_m} )$ is at least $3m+3$. It is clear one may label $C_{3}\circ \overline{K_m}$ so that it is local antimagic, thus $\chi_{la} ( C_{3} \circ \overline{K_m} ) \leq 3m+3$. We are done with the proof.\qed
~~


Now we turn to more general cases. First we obtain the following results for $\chi_{la} ( C_{n} \circ \overline{K_m} )$ whenever $n$ is even.
\begin{lem}
$\chi_{la} ( C_{n} \circ \overline{K_m} )\leq mn+2$ for $n$ even, $n\geq 4$, and $m$ even, $m \geq 4$.
\end{lem}
\proof
Assume the vertex set is $V(C_{n} \circ \overline{K_m}) = \bigcup_{i=1}^{n}\{ d_i^1,d_i^2, \cdots, d_i^m\}\cup\{ u_1,u_2, \cdots, u_n \}$ and the edge set is $E(C_{n} \circ \overline{K_m}) = \bigcup_{i=1}^{n} \{ u_i d_i^j : ~1 \leq i \leq n, 1 \leq j \leq m \}\cup\{ u_i u_{i+1} : ~ 1 \leq i \leq n-1 \} \cup \{u_n u_1 \}$.

With the method used in Lemma~\ref{Cn-K2-n-even}, we first give the labeling $1,2,\cdots, 3n$ to the edges $u_i u_{i+1}, 1 \leq i \leq n-1$, $u_n u_1$, and $u_i d_i^j : ~1 \leq i \leq n, 1 \leq j \leq 2$ for the subgraph $C_{n} \circ \overline{K_2}$. Then the partial vertex sums at $u_i$ is $4n+2$ for odd $i$, and $6n+2$ for even $i$ respectively.

Then secondly we give the edge labels $3n+1, 3n+2, \cdots, mn+n$ to the remaining edges as follows whenever $1\leq i \leq n$ and $3\leq j \leq m$:

\begin{displaymath}
f(u_{i} d_i^{j})= \left\{ \begin{array}{ll}
jn + i, ~~j~ \mbox{odd}.\\
(j+1)n + 1 - i,~~j~ \mbox{even}.
\end{array} \right.
\end{displaymath}
Therefore it can be calculated the vertex sums at $u_i$ is $(\frac{m^2 + 2m}{2})n + \frac{m+2}{2}$ for odd $i$, and $(\frac{m^2 + 2m +4}{2})n + \frac{m+2}{2}$ for even $i$, respectively. Thus we obtain two new colors(vertex sums) other than those over leaves.
We are done with the proof.\qed
~~

Similarly we have the following.
\begin{lem}
$\chi_{la} ( C_{n} \circ \overline{K_m} )\leq mn+2$ for $n$ even, $n\geq 4$ and $m$ odd, $m \geq 3$.
\end{lem}
\proof
Assume the vertex set is $V(C_{n} \circ \overline{K_m}) = \bigcup_{i=1}^{n}\{ d_i^1,d_i^2, \cdots, d_i^m\}\cup\{ u_1,u_2, \cdots, u_n \}$ and the edge set is $E(C_{n} \circ \overline{K_m}) = \bigcup_{i=1}^{n} \{ u_i d_i^j : ~1 \leq i \leq n, 1 \leq j \leq m \}\cup\{ u_i u_{i+1} : ~ 1 \leq i \leq n-1 \} \cup \{u_n u_1 \}$.

Label the edges over the cycle via $f(u_n u_1)=1$ and $u_{i} u_{i+1}, 1\leq i \leq n-1$, as follows:

\begin{displaymath}
f(u_{i} u_{i+1})= \left\{ \begin{array}{ll}
n - \frac{i -1}{2}, i~ \mbox{odd.}\\
2 + \frac{i-2}{2}, i~ \mbox{even.}
\end{array} \right.
\end{displaymath}

Label the edges over $u_i d_i^j$ for $1\leq j \leq 3$ as follows:

\begin{displaymath}
f(u_{i} d_i^{1})= \left\{ \begin{array}{ll}
n+1 + \frac{i -1}{2}, ~~i~ \mbox{odd}.\\
2n+i,~~i~ \mbox{even}.
\end{array} \right.
\end{displaymath}

\begin{displaymath}
f(u_{i} d_i^{2})= \left\{ \begin{array}{ll}
\frac{3n+2}{2} + \frac{i -1}{2}, ~~i~ \mbox{odd}.\\
\frac{7n}{2} - \frac{i-2}{2},~~i~ \mbox{even}.
\end{array} \right.
\end{displaymath}

\begin{displaymath}
f(u_{i} d_i^{3})= \left\{ \begin{array}{ll}
3n - i, ~~i~ \mbox{odd}.\\
4n-1 - \frac{i-2}{2},~~i\neq n,~i~ \mbox{even}.\\
4n, ~~i = n.
\end{array} \right.
\end{displaymath}

Therefore it can be calculated the partial vertex sums at $u_i$ is $\frac{13n+2}{2} + 2$ for odd $i$, and $\frac{21n+2}{2} + 3$ for even $i$, respectively.

Then secondly we give the edge labels $4n+1, 4n+2, \cdots, mn+n$ to the remaining edges as follows whenever $1\leq i \leq n$ and $4\leq j \leq m$:

\begin{displaymath}
f(u_{i} d_i^{1})= \left\{ \begin{array}{ll}
(j+1)n + 1-i, ~~j~ \mbox{odd}.\\
jn+i,~~j~ \mbox{even}.
\end{array} \right.
\end{displaymath}

Then it can be calculated the vertex sums at $u_i$ is $(\frac{m^2 + 2m-2}{2})n + \frac{m+1}{2}$ for odd $i$, and $(\frac{m^2 + 2m +6}{2})n + \frac{m+3}{2}$ for even $i$, respectively.
Thus we obtain two new colors(vertex sums) other than those over leaves.
We are done with the proof.\qed
~~

For odd cycles $C_n, n \geq 5$, we have the following lower bounds.
\begin{lem}\label{C_n_K_m_n_odd_lower_bound}
$\chi_{la} ( C_{n} \circ \overline{K_m} )\geq mn+3$ for $n$ odd and $n < \frac{(m+2)(m+3)}{2(m+1)}$, $m \ge 1$.
\end{lem}
\proof
Assume the vertex set of $C_{n} \circ \overline{K_m}$ is $V = \bigcup_{i=1}^{n}\{ d_i^1,d_i^2, \cdots, d_i^m\}\cup\{ u_1,u_2, \cdots, u_n \}$ and the edge set is $E = \bigcup_{i=1}^{n} \{ u_i d_i^j : ~1 \leq i \leq n, 1 \leq j \leq m \}\cup\{ u_i u_{i+1} : ~ 1 \leq i \leq n-1 \} \cup \{u_n u_1 \}$. Note that there are $mn+n$ edges.

By Lemma~\ref{Cn-Km-lower-bound} we see a lower bound of $\chi_{la} ( C_{n} \circ \overline{K_m} )$ is $mn+2$. Assume that $\chi_{la} ( C_{n} \circ \overline{K_m} )= mn+2$, that is there are two more new colors(vertex sums) other than those over $mn$ leaves in this case. Since $n$ is odd, there must be at least one repeated color, say $w(u_i) = w(d_j^k)$ for some $i\neq j$. On one hand, with smallest edge labels we see that $w(u_i) \geq 1+2+\cdots+(m+2)=\frac{(m+2)(m+3)}{2}$. On the other hand, $w(u_i) = w(d_j^k) \leq mn+n$ which is the largest edge label. Hence we have
$$ mn+n \geq \frac{(m+2)(m+3)}{2} $$
which implies $ n \geq \frac{(m+2)(m+3)}{2(m+1)}$. Therefore we see that if $n < \frac{(m+2)(m+3)}{2(m+1)}$, the local antimagic vertex coloring number  $\chi_{la} ( C_{n} \circ \overline{K_m} )\geq mn+3$. \qed
~~

Note that from the above Theorem~\ref{C3-Km} and the Lemma~\ref{C_n_K_m_n_odd_lower_bound}, we have the following table of lower bounds of local antimagic coloring numbers for smaller $m$'s and $n$'s.\\~~\\

\[
\begin{tabular}{|c|c|c|c|}
\hline
$m$ & $n < \frac{(m+2)(m+3)}{2(m+1)}$ & $n$ & $\chi_{la} ( C_{n} \circ \overline{K_m} )$ \\
\hline
$1$ & $n < 3$ & $.$ & $.$ \\
\hline
$2$ & $n < \frac{20}{6}$ & $3$ & $\geq 2n+3$ \\
\hline
$3$ & $n < \frac{30}{8}$ & $3$ & $\geq 3n+3$ \\
\hline
$4$ & $n < \frac{42}{10}$ & $3$ & $\geq 4n+3$ \\
\hline
$5$ & $n < \frac{56}{12}$ & $3$ & $\geq 5n+3$ \\
\hline
$6$ & $n < \frac{72}{14}$ & $3, 5$ & $\geq 6n+3$ \\
\hline
$7$ & $n < \frac{90}{16}$ & $3, 5$ & $\geq 7n+3$ \\
\hline
$8$ & $n < \frac{110}{18}$ & $3, 5$ & $\geq 8n+3$ \\
\hline
$9$ & $n < \frac{132}{20}$ & $3, 5$ & $\geq 9n+3$ \\
\hline
$10$ & $n < \frac{156}{22}$ & $3, 5, 7$ & $\geq 10n+3$ \\
\hline
$11$ & $n < \frac{182}{24}$ & $3, 5, 7$ & $\geq 11n+3$ \\
\hline
$12$ & $n < \frac{210}{26}$ & $3, 5, 7$ & $\geq 12n+3$ \\
\hline
$13$ & $n < \frac{240}{28}$ & $3, 5, 7$ & $\geq 13n+3$ \\
\hline
$14$ & $n < \frac{272}{30}$ & $3, 5, 7, 9$ & $\geq 14n+3$ \\
\hline
$15$ & $n < \frac{306}{32}$ & $3, 5, 7, 9$ & $\geq 15n+3$ \\
\hline
\end{tabular}
\]

On the other hand we obtain the upper bounds for $C_{n} \circ \overline{K_m}$ as follows, where $n$ is odd.
\begin{lem}
$\chi_{la} ( C_{n} \circ \overline{K_m} )\leq mn+3$ for $n$ odd, $n\geq 5$, and $m$ odd, $m \geq 3$.
\end{lem}
\proof
Assume the vertex set is $V(C_{n} \circ \overline{K_m}) = \bigcup_{i=1}^{n}\{ d_i^1,d_i^2, \cdots, d_i^m\}\cup\{ u_1,u_2, \cdots, u_n \}$ and the edge set is $E(C_{n} \circ \overline{K_m}) = \bigcup_{i=1}^{n} \{ u_i d_i^j : ~1 \leq i \leq n, 1 \leq j \leq m \}\cup\{ u_i u_{i+1} : ~ 1 \leq i \leq n-1 \} \cup \{u_n u_1 \}$.

Label the edges over the subgraph $C_n \circ {K_1}$ with $1, 2, \cdots, 2n$ as in Lemma~\ref{Cn_K1_odd_n}. Therefore it can be calculated the partial vertex sums at $u_i$ is $\frac{5n+7}{2}$ for odd $i\neq n$, $\frac{3n+5}{2}$ for $i=n$, and $\frac{5n+3}{2}$ for even $i$, respectively.

Then secondly we give the edge labels $2n+1, 2n+2, \cdots, mn+n$ to the remaining edges as follows whenever $1\leq i \leq n$ and $2\leq j \leq m$:

\begin{displaymath}
f(u_{i} d_i^{1})= \left\{ \begin{array}{ll}
(j+1)n + 1-i, ~~j~ \mbox{odd}.\\
jn+i,~~j~ \mbox{even}.
\end{array} \right.
\end{displaymath}

Then it can be calculated the vertex sums at $u_i$ is $(\frac{m^2 + 2m+2}{2})n + \frac{m+6}{2}$ for odd $i\neq n$, $(\frac{m^2 + 2m}{2})n + \frac{m+4}{2}$ for $i=n$, and $(\frac{m^2 + 2m +2}{2})n + \frac{m+2}{2}$ for even $i$, respectively.
Thus we obtain three new colors(vertex sums) other than those over leaves.
We are done with the proof.\qed
~~

\begin{lem}
$\chi_{la} ( C_{n} \circ \overline{K_m} )\leq mn+3$ for $n$ odd, $n\geq 5$, and $m$ even, $m \geq 4$.
\end{lem}
\proof
Assume the vertex set is $V(C_{n} \circ \overline{K_m}) = \bigcup_{i=1}^{n}\{ d_i^1,d_i^2, \cdots, d_i^m\}\cup\{ u_1,u_2, \cdots, u_n \}$ and the edge set is $E(C_{n} \circ \overline{K_m}) = \bigcup_{i=1}^{n} \{ u_i d_i^j : ~1 \leq i \leq n, 1 \leq j \leq m \}\cup\{ u_i u_{i+1} : ~ 1 \leq i \leq n-1 \} \cup \{u_n u_1 \}$.

Label the edges over the subgraph $C_n \circ \overline{K_2}$ with $1, 2, \cdots, 3n$ as in Lemma~\ref{Cn-K2-odd}. Therefore it can be calculated the partial vertex sums at $u_i$ is $5n+6$ for odd $i\neq n$, $\frac{3n+13}{2}$ for $i=n$, and $5n+7$ for even $i$, respectively.

Then secondly we give the edge labels $3n+1, 3n+2, \cdots, mn+n$ to the remaining edges as follows whenever $1\leq i \leq n$ and $3\leq j \leq m$:

\begin{displaymath}
f(u_{i} d_i^{1})= \left\{ \begin{array}{ll}
jn + i, ~~j~ \mbox{odd}.\\
(j+1)n+1-i,~~j~ \mbox{even}.
\end{array} \right.
\end{displaymath}

Then it can be calculated the vertex sums at $u_i$ is $(\frac{m^2 + 2m+2}{2})n + \frac{m+10}{2}$ for odd $i\neq n$, $(\frac{m^2 + 2m-5}{2})n + \frac{m+11}{2}$ for $i=n$, and $(\frac{m^2 + 2m +2}{2})n + \frac{m+12}{2}$ for even $i$, respectively.
Thus we obtain three new colors(vertex sums) other than those over leaves.
We are done with the proof.\qed
~~

To summarize with previous observations, we have the following results.
\begin{thm}
For each odd cycle $C_n, n\geq 5$, the local antimagic chromatic number $\chi_{la} ( C_{n} \circ \overline{K_m} ) = mn+3$ except finitely many $m$'s.
\end{thm}

As seen in Lemma~\ref{C_n_K_m_n_odd_lower_bound} and the associated table, we may precisely have $\chi_{la} ( C_{5} \circ \overline{K_m} ) = mn+3$ for $m \geq 6$, also $mn+2 \leq \chi_{la} ( C_{5} \circ \overline{K_m} ) \leq mn+3$ for $3\leq m \leq 5$. The next one is $\chi_{la} ( C_{7} \circ \overline{K_m} ) = mn+3$ for $m \geq 10$, also $mn+2 \leq \chi_{la} ( C_{7} \circ \overline{K_m} ) \leq mn+3$ for $3\leq m \leq 9$ and so on.

\section{Corona Product of Complete Graphs with Null Graphs}
\subsection{The Chromatic Number of $K_{n} \circ K_1$}

\begin{thm}
$\chi_{la} ( K_{n} \circ K_1 )= 2n-1$ for $n\ge 3$.
\end{thm}

\proof

\begin{figure}[h]
\centering
\includegraphics[width=14cm]{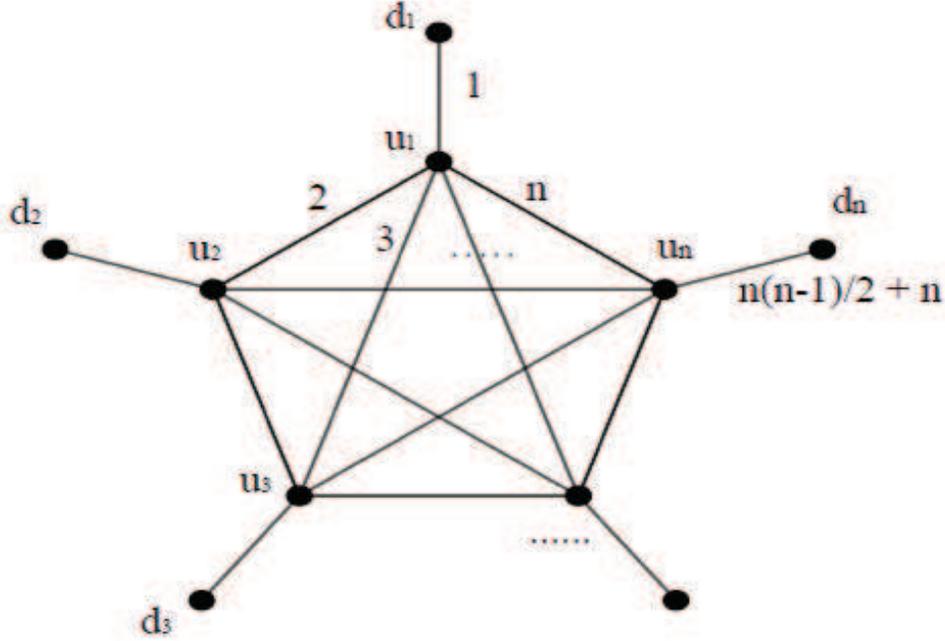}
\caption{Corona product $K_{n} \circ K_1$}
\label{kn-k1}
\end{figure}

The corona product graph $K_{n} \circ K_1$ has $2n$ vertices and $\binom{n}{2}+n$ edges. Suppose the vertex set is $\{ u_1, \cdots, u_n \} \cup \{ d_1, \cdots, d_n \}$, where $d_i$ is of degree one and $u_i$ is of degree $n$, for $i=1,\cdots,n$. Note that the pendent edges are $u_1 d_1, u_2 d_2, \cdots, u_n d_n$. See Figure~\ref{kn-k1} for reference.

First we see $2n-1$ is the lower bound of $\chi_{la} ( K_{n} \circ K_1 )$.
For vertices $u_i$ the vertex sum(color) $w(u_i)$ can be estimated as follows for each $i$:
$$ w(u_i) \geq 1+2+\cdots+(n-1)+n = \binom{n+1}{2} = \binom{n}{2}+n$$
Note that $\binom{n}{2}+n$ is the largest edge label used. Thus it implies at most one vertex $u_i$ could have the same repeated color with $w(d_j)$ for some $j$, where $w(u_i)=1+2+\cdots+(n-1)+n = \binom{n}{2}+n = w(d_j)$. Therefore $\chi_{la} ( K_{n} \circ K_1 ) \geq 2n-1$ for $n\ge 3$, which gives a lower bound.

On the other hand, one can see $2n-1$ is also the upper bound of $\chi_{la} ( K_{n} \circ K_1 )$. Suppose without loss of generality at the vertex $u_1$ the incident edges are labeled $1, 2, \cdots, n$ so that $u_1 d_1$ is labeled as 1(see Figure~\ref{kn-k1}).
Let us hold the largest $n-1$ edge labels $\binom{n}{2}+2, \binom{n}{2}+3, \cdots, \binom{n}{2}+n-1, \binom{n}{2}+n$ for later use. Then arrange arbitrarily $n+1, n+2, \cdots, \binom{n}{2}+1$ on the edges within the $K_n$ part, and sort the partial vertex sums(without loss of generality again) as $\widetilde{w}(u_2) \leq \widetilde{w}(u_3) \leq \cdots \leq \widetilde{w}(u_n)$, where $\widetilde{w}(u_i)$ is the partial vertex sum of incident edge labels without the pendent edge label at $u_i$. Now put the largest $n-1$ edge labels $\binom{n}{2}+2, \binom{n}{2}+3, \cdots, \binom{n}{2}+n-1, \binom{n}{2}+n$ over the pendent edges $u_2 d_2, u_3 d_3, \cdots, u_n d_n$ respectively. Therefore we obtain an edge labeling so that $w(d_1)=1$,  $w(u_1)=\binom{n}{2}+n=w(d_n)$, and $w(u_2) < \cdots < w(u_n)$. Note that the remaining colors $w(d_2), w(d_3), \cdots, w(d_n)$ are respectively $\binom{n}{2}+2, \binom{n}{2}+3, \cdots, \binom{n}{2}+n-1, \binom{n}{2}+n$. Also one sees that $w(u_i) > \binom{n}{2}+n$ for each $i \geq 2$. That is we have the colors $1=w(d_1)<w(d_2)< \cdots < w(d_n)=\binom{n}{2}+n = w(u_1) < w(u_2) < \cdots < w(u_n)$, which gives the upper bound $\chi_{la} ( K_{n} \circ K_1 ) \leq 2n-1$ and we are done with the proof.
\qed
~~

\subsection{The Chromatic Number of $K_{n} \circ \overline{K_m}$}

Note that the corona product graph $K_{n} \circ \overline{K_m}$ has $mn+n$ vertices and $\binom{n}{2}+mn$ edges. Suppose the vertex set is $\bigcup_{i=1}^{n} \{ d_{i1}, \cdots, d_{im} \}\cup \{ u_1, \cdots, u_n \} $, where $d_{ij}$ are leaves of degree one and $u_i$ are of degree $m+(n-1)$, for $i=1,\cdots,n$ and $j=1,\cdots,m$. More generally we have the following.
\begin{thm}
$\chi_{la} ( K_{n} \circ \overline{K_m} )= mn+n$ for $m\ge 2, n\ge 3$.
\end{thm}

\proof

We claim that for each $u_i$, the vertex sum(color) $w(u_i)$ will not be repeated with the same vertex sum(color) as any leaf $w(d_{jk})$, for any $1\leq j \leq n, 1\leq k \leq m$. Therefore this implies $\chi_{la} ( K_{n} \circ \overline{K_m} )= mn+n$ immediately since the colors $w(u_i)$ are all distinct and are new relative to those leaf-colors.

To show the claim, suppose on the contrary. We see first that using the smallest edge labels $w(u_i) \geq 1+2+\cdots+(m+n-1)$. On the other hand, since we assume for some $u_i$, its color $w(u_i)$ is repeated with some $w(d_{jk})$, therefore $w(u_i) \leq \binom{n}{2}+mn$ which is the largest edge label used. Combining the above we have
$$1+2+\cdots+(m+n-1) \leq \binom{n}{2}+mn.$$
After simplification one has $m(m-1) \leq 0$, which is a contradiction since $m \geq 2$. We are done with the proof.\qed
~~

\section{Concluding Remarks}

In this article among others we determine completely the local antimagic chromatic number $\chi_{la}(G\circ \overline{K_m})$ for the corona product of a graph $G$ with the null graph $\overline{K_m}$ on $m\geq 1$ vertices, when $G$ is a path $P_n$, a cycle $C_n$, and a complete graph $K_n$.

Just recently G.C. Lau, H.K. Ng, W.C. Shiu\cite{LNS1, LNS2, LNS3} and S. Shaebani\cite{S} respectively studied the concept local antimagic chromatic numbers for various classes of graphs for join products. More variants of the concept extended to edge coloring and total coloring can be explored.

\end{document}